\documentclass[11pt]{amsart}
\usepackage{geometry}                
\usepackage{graphicx}
\usepackage{amssymb}
\usepackage{epstopdf}
\usepackage{longtable}
\usepackage{psfrag}
\usepackage{overpic}
\usepackage{color}
\usepackage{caption2}
\usepackage{subfigure}

\title{Two elementary formulae and some complicated properties for Mertens function}

\author{Rong Qiang Wei}
\address{College of Earth Sciences, University of Chinese Academy of Sciences, Beijing, PRC, 100049}
\email{wrq1973@ucas.ac.cn}
\date{}

\begin{document}
\maketitle
\begin{abstract}
     Two elementary formulae for Mertens function $M(n)$ are obtained. With these formulae, $M(n)$ can be calculated directly and simply, which can be easily implemented by computer. $M (1) \sim M (2 \times 10^7) $ are calculated one by one. Based on these $2\times 10^7$ samples, some of the complicated properties for Mertens function $M(n)$,  its 16479 zeros, the 10043 local maximum/minimum between two neighbor zeros,  and the relation with the cumulative sum of the squarefree integers, are understood numerically and empirically.
  
\end{abstract}

{\hspace{2.2em}\small Keywords:}

{\hspace{2.2em}\tiny Mertens function Elementary formula Zeros Local maximum/minimum squarefree integers}

\section{Introduction}

The Mertens function $M(n)$ is important in number theory, especially the calculation of this function when its argument is large. $M(n)$ is defined as the cumulative sum of the M$\ddot{\mathrm{o}}$bius function $\mu(k)$ for all positive integers $n$,

\begin{equation}\label{eq1}
M(n)=\sum_{k=1}^{n}\mu(k)
\end{equation}

where the M$\ddot{\mathrm{o}}$bius function $\mu(k)$ is defined for a positive integer $k$ by

\begin{equation}\label{eq2}
\mu(k)  = \left\{ {\begin{array}{*{20}{c}}
1\\
0\\
{{{( - 1)}^m}}
\end{array}\begin{array}{*{20}{c}}
{k = 1{\mbox{\hspace{14em} }}}\\
{{\mbox{if \ }} k {\mbox{\ is divisible by a prime square\hspace{2em}}}}\\
{{\mbox{if\  }} k {\mbox{\  is the product of }}m{\mbox{ distinct primes}}}
\end{array}} \right.
\end{equation}

Sometimes the above definition can be extended to real numbers as follows:

\begin{equation}\label{eq3}
M(x)=\sum_{1\le k \le x}\mu(k)
\end{equation}

Furthermore, Mertens function has other representations. They are mainly shown in formula (\ref{eq4}-\ref{eq7}):

\begin{equation}\label{eq4}
M(x)=\frac{1}{2\pi i}\int_{c-i\infty}^{c+i\infty}\frac{x^s}{s\zeta (s)}\mbox{d}s
\end{equation}

where $\zeta (s)$ is the the Riemann zeta function, $s$ is complex, and $c>1$

\begin{equation}\label{eq5}
M(n)=\sum_{a\in F_n}\exp(2\pi i a)
\end{equation}

where $F_n$ is the Farey sequence of order $n$.

\begin{equation}\label{eq6}
M(n)=1-\sum_{a=2}^{a\le n}1+\sum_{a=2}^{a\le\frac{n}{b}}\sum_{b=2}^{b\le\frac{n}{a}}1-\sum_{a=2}^{a\le\frac{n}{bc}}\sum_{b=2}^{b\le\frac{n}{ac}}\sum_{c=2}^{b\le\frac{n}{ab}}1+\sum_{a=2}^{a\le\frac{n}{bcd}}\sum_{b=2}^{b\le\frac{n}{acd}}\sum_{c=2}^{b\le\frac{n}{abd}}\sum_{d=2}^{b\le\frac{n}{abc}}1\dots
\end{equation}

Formula (\ref{eq6}) shows $M(n)$ is the sum of the number of points under n-dimensional hyperboloids.

\begin{equation}\label{eq7}
M(n)=\mbox{det}R_{n\times n}
\end{equation}

where $R$ is the Redheffer matrix. $R = \{ r_{i j} \}$ is defined by $r_{i j} = 1$ if $j =1$ or
$i \vert j $, and $r_{i j} = 0$ otherwise.

However, none of the representations above results in practical algorithms for calculating the Mertens function. At present, many algorithms are based or partially based on sieving similar to those used in prime counting (eg., Kotnik and van de Lune, 2004; Kuznetsov, 2011; Hurst, 2016).  With the sieve algorithm, $M(x)$ has been computed for all $x\le 10^{22}$(Kuznetsov, 2011). For the isolated values of Mertens function $M(n)$, $M(2^n)$ has been computed for all positive integers $n\le 73$ with combinatorial algorithm (partially sieving) (Hurst, 2016).

On the other hand, there are also some (recursive) formulae for Mertens function $M(n)$ in the mathematical literature by which the practical algorithms for the $M(n)$ can be obtained (eg., Neubauer, 1963; Dress, 1993; Benito and Varona 2008; and the references therein).  For example, Benito and Varona (2008) present a two-parametric family of recursive formula as follows,

\begin{equation}\label{eq7_add}
\begin{array}{*{20}{l}}
{2M(n) + 3}&{ = \sum\limits_{k = 1}^{\left\lfloor {\frac{n}{{r + 1}}} \right\rfloor } {g(n,k)\mu (k)} }\\
{}&{ + \sum\limits_{b = 0}^s {[M(\frac{n}{{3 + 6b}}) - 2M(\frac{n}{{5 + 6b}}) + 3M(\frac{n}{{6 + 6b}}) - 2M(\frac{n}{{7 + 6b}})]} }\\
{}&{ + \sum\limits_{a = 6s + 9}^r {h(a)[} M(\frac{n}{a}) - M(\frac{n}{{a + 1}})]}
\end{array}
\end{equation}

where $n$, $r$, and $s$ are three integers such that $s\ge 0$ and $6s+9\le r \le n-1$. $g(n,k)=3\lfloor\frac{n}{3k}\rfloor -2\lfloor\frac{n}{2k}-\frac{1}{2}\rfloor$. $h(a)=g(n,k)$ for $\frac{n}{a+1}<k\le\frac{n}{a}$.

Here we introduce two new elementary formulae to calculate the $M(n)$ which are not based on sieving.  We calculate Mertens function from $M (1) $ to $M(2 \times 10^7)$ with these formulae one by one,  and study some properties of the $M(n)$ numerically and empirically.

\section{Two elementary formulae for Mertens function}\label{sec2}

In Wei (2016), a definite recursive relation for M$\ddot{\mathrm{o}}$bius function is
introduced by two simple ways. One is from M$\ddot{\mathrm{o}}$bius transform, and the
other is from the submatrix of the Redheffer Matrix.

The recursive relation for M$\ddot{\mathrm{o}}$bius function $\mu(k)$ is,

\begin{equation}\label{eq8}
\mu (k) =  - \sum\limits_{m = 1}^{k - 1} {{l_{km}}} \mu (m),k = 2,3, \cdots ;{l_{km}} = \left\{ {\begin{array}{*{20}{c}}
1&{m|k}\\
0&{{\rm{else}}}
\end{array}} \right.
\end{equation}

and $\mu(1)=1$.

From formula (\ref{eq8}), two elementary formulae for Mertens function $M(n)$ can be obtained as follows,

\begin{equation}\label{eq9}
M(n) = \sum_{k=1}^{n}\mu(k)= 1+\sum_{k=2}^{n}[- \sum\limits_{m = 1}^{k - 1} {{l_{km}}} \mu (m)],k = 2,3, \cdots ;{l_{km}} = \left\{ {\begin{array}{*{20}{c}}
1&{m|k}\\
0&{{\rm{else}}}
\end{array}} \right.
\end{equation}

\begin{equation}\label{eq9_add}
M(n) = M(n-1)- \sum\limits_{m = 1}^{n - 1} {{l_{nm}}} \mu (m),n = 2,3, \cdots ;{l_{nm}} = \left\{ {\begin{array}{*{20}{c}}
1&{m|n}\\
0&{{\rm{else}}}
\end{array}} \right.
\end{equation}

where $M(1)$=1.

With formula (\ref{eq9}) or (\ref{eq9_add}), $M(n)$ can be calculated directly and the most complex operation is only the Mod. We calculated the Mertens function from $M (1)$ to $M (2 \times 10^7)$ one by one. Some values of the $M(n)$ are listed in Table \ref{tb1}.

\begin{table}[htdp]
\small
\caption{Some values of the Mertens function $M(n)$}
\begin{center}
\begin{tabular}{cccccccccccc}
\hline
$n$   &$M(n)$& $n$ & $M(n)$& $n$ & $M(n)$ &$n$ &$M(n)$&$n$ &$M(n)$&$n$ &$M(n)$\\
\hline
 10 & -1	& 20           &-3	 & 30           &  -3	& 40           &0	  & 50           & -3	 & 60           & -1    \\   
$10^2$&1	&$2\times 10^2$&-8	 &$3\times 10^2$&  -5	&$4\times 10^2$&1	  &$5\times 10^2$& -6	 &$6\times 10^2$&  4    \\      
$10^3$&2	&$2\times 10^3$& 5	 &$3\times 10^3$&  -6	&$4\times 10^3$&-9  &$5\times 10^3$&  2	 &$6\times 10^3$&  0    \\        
$10^4$&-23&$2\times 10^4$& 26	 &$3\times 10^4$&  18	&$4\times 10^4$&-10	&$5\times 10^4$&  23 &$6\times 10^4$&  -83  \\  
$10^5$&-48&$2\times 10^5$& -1	 &$3\times 10^5$&  220&$4\times 10^5$&11	&$5\times 10^5$& -6	 &$6\times 10^5$&  -230 \\  
$10^6$&212&$2\times 10^6$& -247&$3\times 10^6$&	 107&$4\times 10^6$&192	&$5\times 10^6$& -709&$6\times 10^6$&	 257  \\
$10^7$&1037&$2\times 10^7$& -953&$3.5\times 10^6$&	 -138&$4.5\times 10^6$&173	&$5.5\times 10^6$& -513&$6.5\times 10^6$&	 867  \\
\hline
\end{tabular}
\end{center}
\label{tb1}
\end{table}

\section{Some complicated properties for Mertens function $M(n)$}

\subsection{Variation of Mertens function $M(n)$ with $n$}

\begin{figure}[htb]
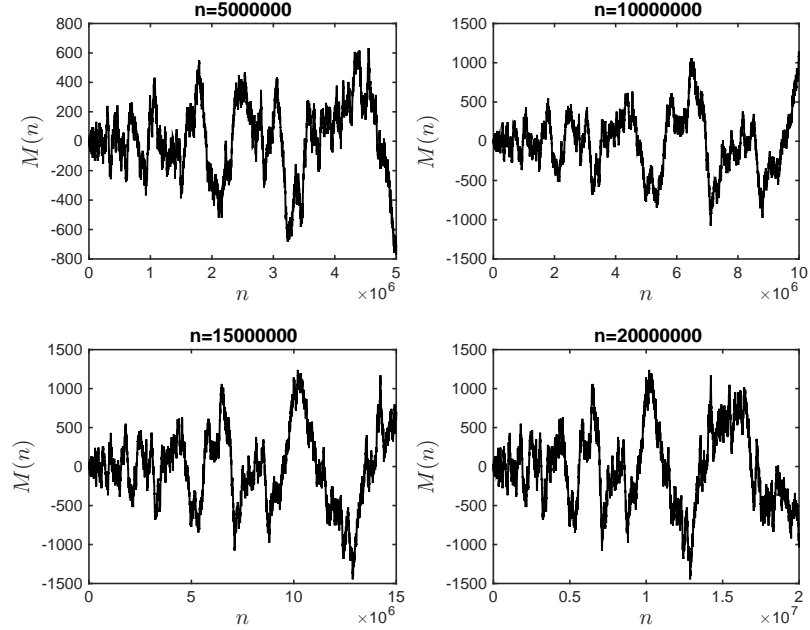

\setlength{\belowcaptionskip}{0pt}
\centering
\begin{overpic}[scale=0.6]{Fig1_m.eps}
\end{overpic}
\renewcommand{\figurename}{Fig.}
\caption{Mertens function $M(n)$ to $n=5\times 10^5, 1\times 10^6, 1.5\times 10^7, 2\times 10^7$.}
\label{fig1}
\end{figure}

Figure \ref{fig1} shows the Mertens function $M(n)$ to $n=5\times 10^5$, $n=1\times 10^6$, $n=1.5\times 10^7$, and $n=2\times 10^7$, respectively. It can be found the distribution of Mertens function $M(n)$ is complicated. It oscillates up and down with increasing amplitude over $n$, but grows slowly in the positive or negative directions until it increases to a certain peak value.

\begin{figure}[htb]
\setlength{\belowcaptionskip}{0pt}
\centering
\begin{overpic}[scale=0.6]{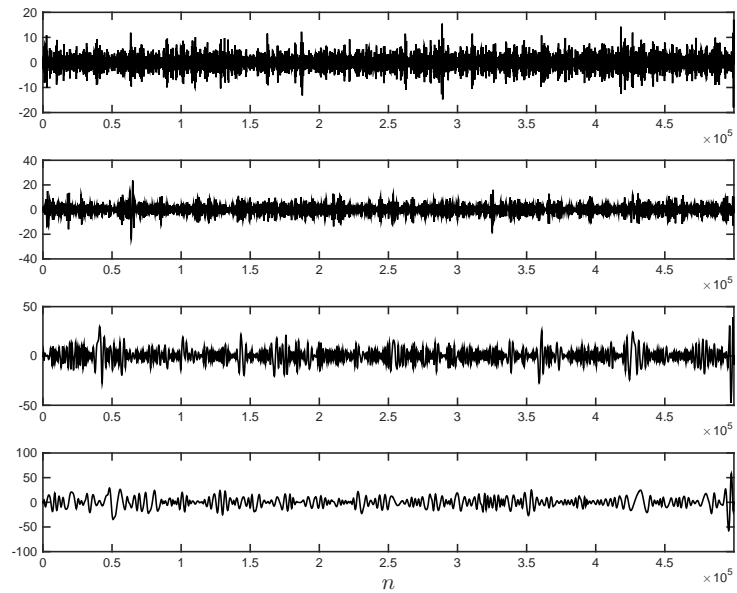}
\end{overpic}
\renewcommand{\figurename}{Fig.}
\caption{The 9th-12th (from the top to the bottom) empirical modes of the sequence composed of $M(1)-M(5\times 10^5)$}
\label{fig2}
\end{figure}

\begin{figure}[htb]
\setlength{\belowcaptionskip}{0pt}
\centering
\begin{overpic}[scale=0.6,bb=50 126 517 402]{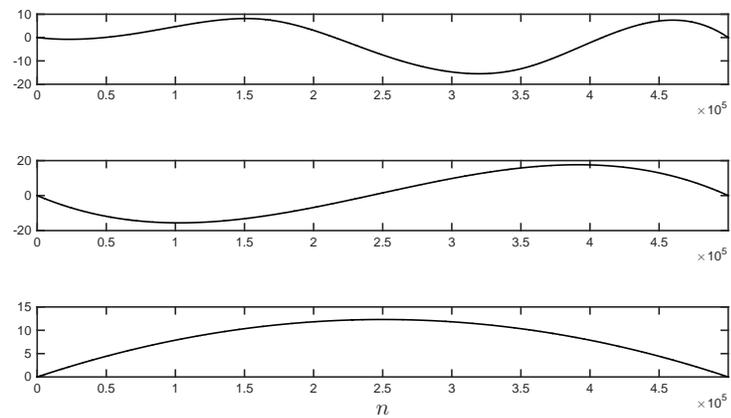}
\end{overpic}
\renewcommand{\figurename}{Fig.}
\caption{The 17th-19th (from the top to the bottom) empirical modes of the sequence composed of $M(1)-M(5\times 10^5)$}
\label{fig3}
\end{figure}

To investigate further the complicated properties of the $M(n)$, a sequence composed of $M(1)-M(5\times 10^5)$ is analyzed by the method of empirical mode decomposition (EMD) (Tan, 2016).  This sequence is decomposed into a sum of 19 empirical modes.  Some models are shown in Figure \ref{fig2}-\ref{fig3}. It can be seen that the empirical modes are still complicated except 16th-19th modes. The fundamental mode (the last mode in Figure \ref{fig3}) is a simple parabola going downward. These empirical modes show further that the Mertens function $M(n)$ has complicated behaviors.

%

\begin{figure}[htb]
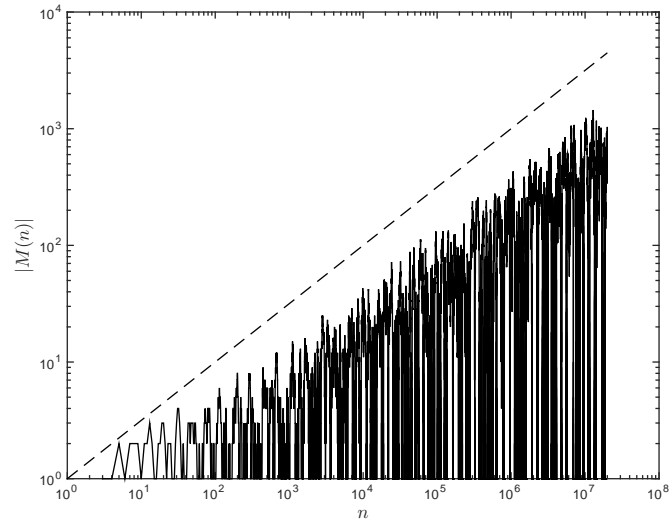

\setlength{\belowcaptionskip}{0pt}
\centering
\begin{overpic}[scale=0.5]{Fig4_m.eps}
\end{overpic}
\renewcommand{\figurename}{Fig.}
\caption{The absolute value of the Mertens function $M(n)$ to $n=2\times 10^7$. The dash line is $\sqrt{n}$.}
\label{fig4}
\end{figure}

It can be found that there is a different variation of the Mertens function $M(n)$ with $n$ in a log-log space. Figure \ref{fig4} shows the absolute value of the Mertens function $\vert M(n)\vert$ to $n=2\times 10^7$ in such a space. It can be found $\vert M(n)\vert$ still oscillates but increases generally with $n$. It seems that the outermost values of $\vert M(n)\vert$ increase almost linearly with $n$, but they are less than those from $\sqrt{n}$.

\begin{figure}[htb]
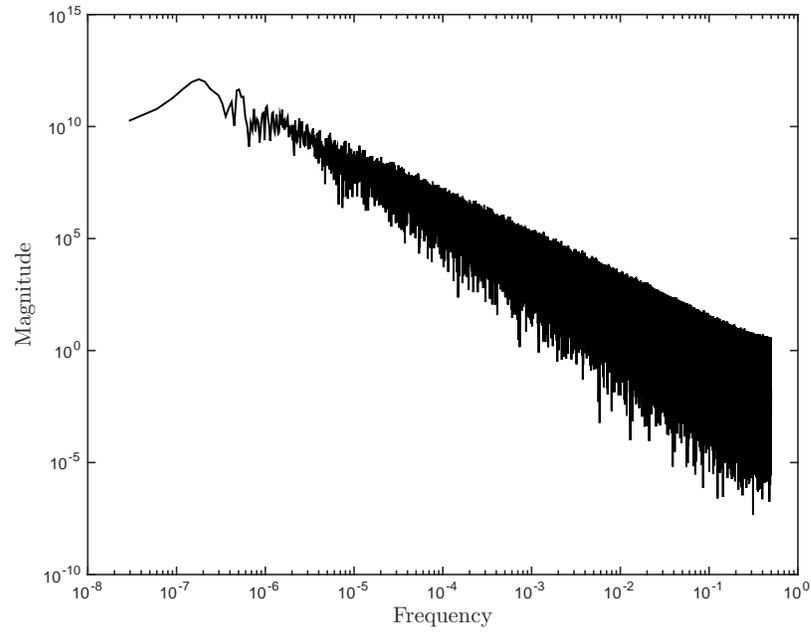

\setlength{\belowcaptionskip}{0pt}
\centering
\begin{overpic}[scale=0.6]{Fig4_add.eps}
\end{overpic}
\renewcommand{\figurename}{Fig.}
\caption{The power spectral density (PSD) for $M(n)$ sequence from $M(1)$ to $M(2\times 10^7)$.}
\label{fig4_add}
\end{figure}

It seems from Figure \ref{fig1} and \ref{fig4} to the $M(n)$ has some fractal properties.  We check this by estimating $M(n)$'s power spectral density (PSD). We calculate the PSD for $M(n)$ sequence from $M(1)$ to $M(2\times 10^7)$ by taking $n$ as time. The result is shown in Figure \ref{fig4_add}.  It can be found that the logarithmic PSD of $M(n)$ has a linear decreasing trend with increasing logarithmic frequency, which does indicate that $M(n)$ has some fractal properties.  However, the M$\ddot{\mathrm{o}}$bius function $\mu(n)$, which is taken as an independent random sequence in Wei (2016), has no such properties. Its cumulative sum, $M(n)$ reduces partly the randomness. 

\subsection{Zeros of the Mertens function $M(n)$}

\begin{figure}[htb]
\setlength{\belowcaptionskip}{0pt}
\centering
\begin{overpic}[scale=0.6, bb=23 44 509 428]{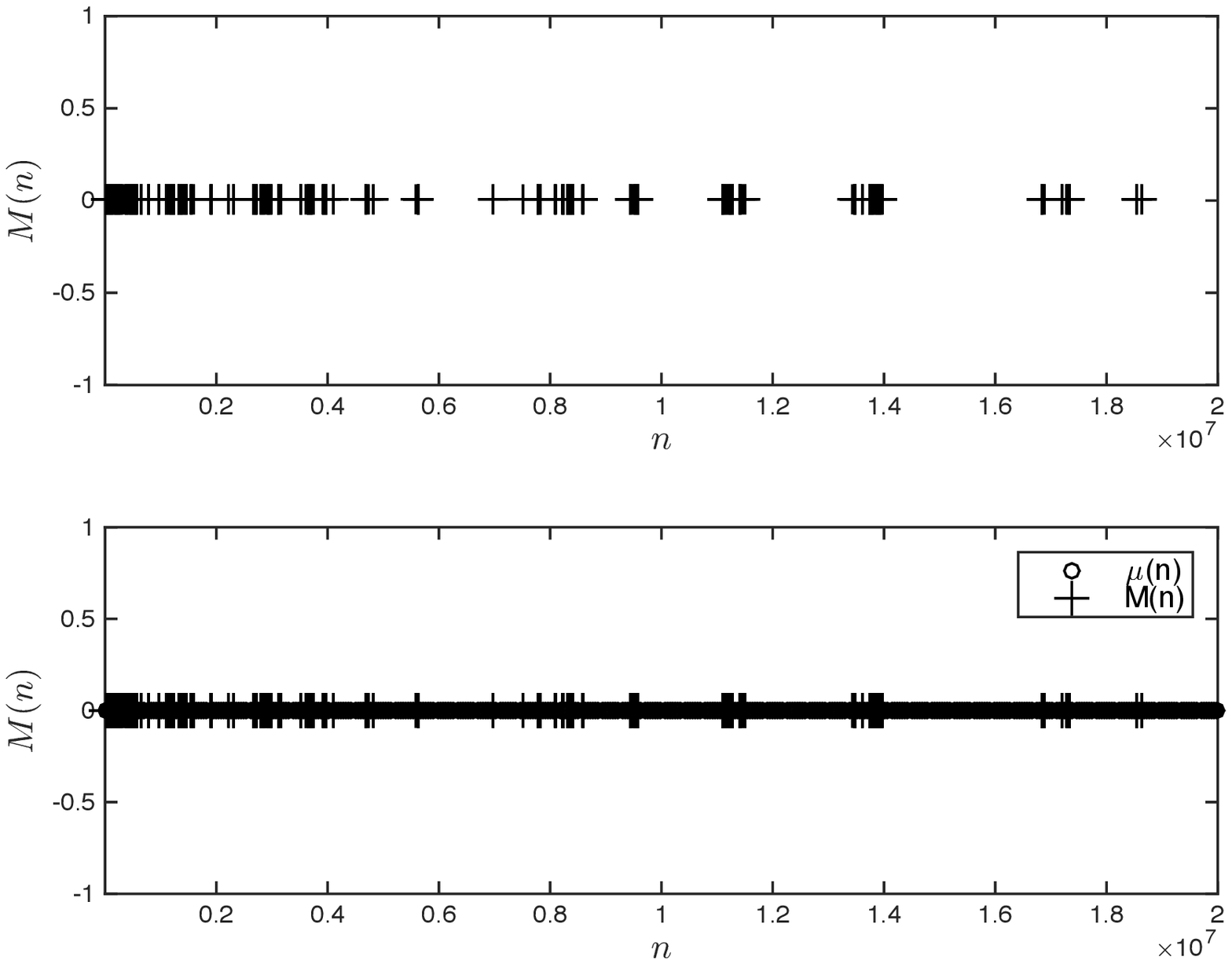}
\end{overpic}
\renewcommand{\figurename}{Fig.}
\caption{Distribution of the zeros of the Mertens function $M(n)$ when $n\in [1, 2\times 10^7]$ (Top);  Comparison of the zeros of the Mertens function $M(n)$ with those of the M$\ddot{\mathrm{o}}$bius function $\mu(n)$(Bottom). }
\label{fig5}
\end{figure}

It is interesting to investigate the distribution of the zeros of the Mertens function $M(n)$. Our calculation shows that there are 16479 zeros among these $2\times 10^7$ values of the $M(n)$, which are shown in the top subfigure of the Figure \ref{fig5}. It can be seen that the zeros are uneven in our computational domain.  Zeros are denser from 1 to $4\times 10^6$, and at $\sim 8\times 10^6$.  The bottom subfigure of the Figure \ref{fig5} shows the comparison of the zeros of the Mertens function $M(n)$ and those of the M$\ddot{\mathrm{o}}$bius function $\mu(n)$ in our computational domain. One can see that the (7841425) zeros of the $\mu(n)$ are much denser than those of the $M(n)$.

\subsection{The  local maximum/minimum of the Mertens function $M(n)$ between two neighbor zeros}

\begin{figure}[htb]
\setlength{\belowcaptionskip}{0pt}
\centering
\begin{overpic}[scale=0.6, bb=22 50 509 429]{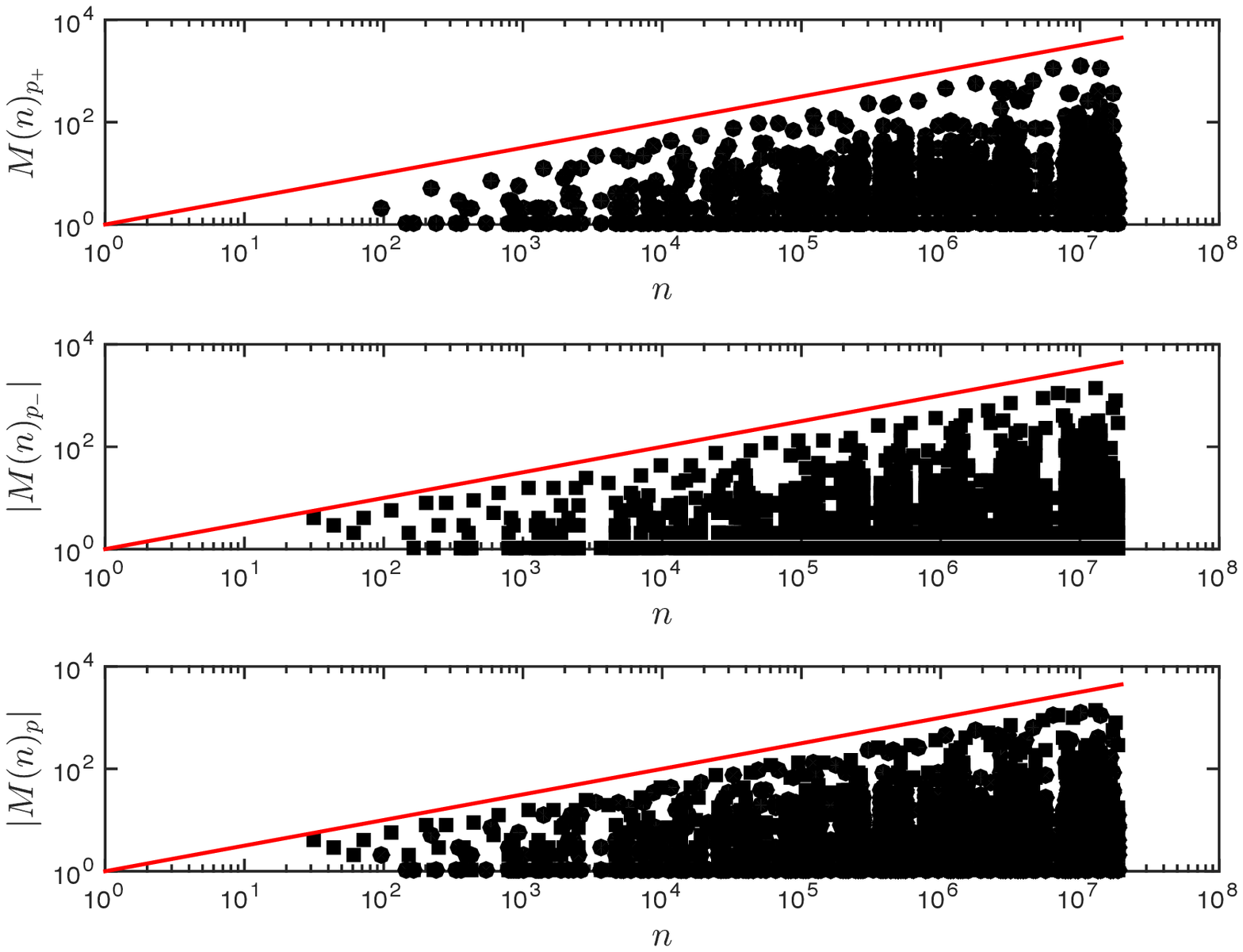}
\end{overpic}
\renewcommand{\figurename}{Fig.}
\caption{Variation of the local maximum/minimum of the $M(1)-M(2\times 10^7)$ with $n$ between two neighbor zeros. The red straight line in each subfigure is $\sqrt{n}$. Top: 5040 maximum values of $M(1)-M(2\times 10^7)$; Middle: 5003 minimum values; Bottom: All the maximum/minimum (solid circle: maximum; solid square: minimum). In the middle and bottom subfigures, the absolute of the minimums are used.}
\label{fig6}
\end{figure}

Another important property for the Mertens function $M(n)$ is the local maximum/minimum of the $M(n)$ between two neighbor zeros. Here the "local" means that the sequence of $M(n)$ must have three values at least including the neighbor zeros.  Figure \ref{fig6} shows the variation of such 10043 local maximum/minimum of $\vert M(n)\vert$ (5040 positive values and 5003 minus values) from $M(1)$ to $M(2\times 10^7)$ with $n$. One can see that the logarithms of the these local maximum/minimum have the similar properties to $\log(\vert M(n)\vert)$ vs. $\log(n)$ in Figure \ref{fig4}. It should be pointed out that the local maximum/minimum here is counted only once when duplicate or more $n$ have the same local maximum/minimum between two neighbor zeros. 

The maximum of these $2\times 10^7$ $M(n)$s is 1240 when $n=10195458$, 10195467, 10195468, 10195522; And the 
minimum is $-1447$ when  $n=12875814$, 12875815, 12875816, 12875818.

\section{Discussions}

\subsection{The calculation of $M(n) $ with the relation (\ref{eq9}) or (\ref{eq9_add})}

In theory, we can calculate $M(n)$ for any large $n$ with the relation (\ref{eq9}) or (\ref{eq9_add}) obtained in section \ref{sec2}. The algorithm is not complicated and can be implemented easily, even by hand. The most operation is only Mod. 
However, in order to calculate $M(n)$ with the relation (\ref{eq9}), it demands $\mu (2), \mu (3),
\ldots  \mu (n - 1)$ firstly, or $M(1)$, $M(2)$, $\dots$, $M(n-1)$ if the relation (\ref{eq9_add}) is used.  Above all, we have to recalculate $l_{km}$ or $l_{nm}$ firstly when $k$ or $n$ changes. It will take a lot of calculation time, especially when $n$ is large. It took about 78326s to get $\mu(5500000)$.  It can be also found that the relation here is not efficient for the isolated $M(n)$. 

In this paper, we only calculate the values of $M(n)$ from $M(1)$ to M($2 \times 10^7)$ because of the limitation of our desktop computer and computing time.  To obtain more numerical results of $M(n)$ with large $n$, both the faster and/or optimization algorithm for the relation (\ref{eq9}) or (\ref{eq9_add}) here, or other professional and efficient algorithms, are required.  

\subsection{Upper bound of  $M(n)$ sequence}

In Wei (2016), the upper bound of  $M(n)$ sequence is discussed based on the assumption that $\mu (n)$ is an independent random sequence, because of the numerical consistency between empirical statistical quantities for only $2\times 10^7$ $\mu(n)$ and those from number theory. The following inequality (\ref{eq10}) for $M(n)$ holds with a probability of $1 - \alpha$, 

\begin{equation}\label{eq10}
M(n) \le  \sqrt{\frac{6}{\pi^2}} K_{_{\frac{\alpha}{2}}}\sqrt{n}
\end{equation}

where,

\begin{equation}
\int\limits_{ - {K_{\alpha /2}}}^{{K_{\alpha /2}}} {\frac{1}{{\sqrt {2\pi } }}} \exp ( - \frac{{{t^2}}}{2}){\rm{d}}t = 1 - \alpha 
\end{equation}

or, the following inequality (\ref{eq10_1}) holds with a probability $p>1-\alpha$ 

\begin{equation}\label{eq10_1}
M(n) \le  \frac{\sqrt {6/{\pi ^2}}}{\sqrt{\alpha}} \sqrt{n}
\end{equation}

\begin{figure}[htb]
\setlength{\belowcaptionskip}{0pt}
\centering
\begin{overpic}[scale=0.4, ]{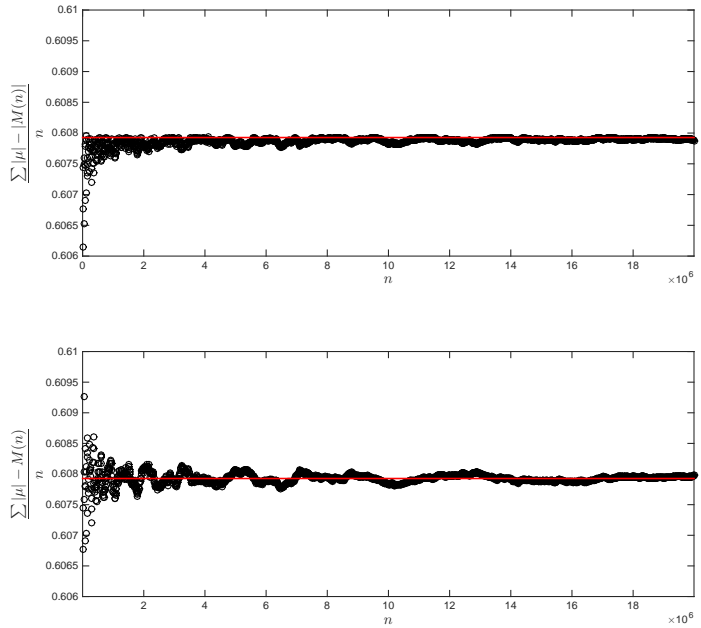}
\end{overpic}
\renewcommand{\figurename}{Fig.}
\caption{Top: Variation of $\frac{\sum\vert\mu(n)\vert}{n} -\frac{\vert M(n)\vert}{n}$ with $n$ from $n=1$ to $n=2\times 10^7$. Bottom: Variation of $\frac{\sum\vert\mu(n)\vert}{n}-\frac{\sum\mu(n)}{n}$ with $n$ from $n=1$ to $n=2\times 10^7$. The red line in each subfigure is $\frac{6}{\pi^2}$.}
\label{fig7}
\end{figure}

Here without taking $\mu (n)$ as an independent random sequence, we conjecture that the upper bound of the $M(n)$ should have a similar formula to (\ref{eq10}) or (\ref{eq10_1}) but the coefficient before $\sqrt{n}$ will be very large from some facts as the follows: 

{\bf Fact 1}.  According to Hardy and Wright (2008),  we have,

\begin{equation}\label{eq11}
\sum\vert\mu(n)\vert=\frac{6}{\pi^2}n+O(n^{1/2})
\end{equation}

\begin{equation}\label{eq12}
M(n)=\sum\mu(n)=o(n)
\end{equation}

For a large $n$, from (\ref{eq11}) we have,

\begin{equation}\label{eq13}
\arrowvert \sum\vert\mu(n)\vert - \frac{6}{\pi^2}n\arrowvert =\sum\vert\mu(n)\vert - \frac{6}{\pi^2}n \leq An^\frac{1}{2}
\end{equation}

where $A$ is a constant.

Then,

\begin{equation}\label{eq15}
\frac{\sum\vert\mu(n)\vert}{n} - \frac{6}{\pi^2}\leq An^{-\frac{1}{2}}
\end{equation}

According to (\ref{eq12}), $\vert M(n)\vert=o(n)$. We have,

\begin{equation}\label{eq16}
\frac{\sum\vert\mu(n)\vert}{n} - \frac{6}{\pi^2}-\frac{\vert M(n)\vert}{n}\leq A_1n^{-\frac{1}{2}}
\end{equation}

where $A_1<A$ is another constant for a large $n$. Further, 

\begin{equation}\label{eq17}
\sum\vert\mu(n)\vert - \frac{6}{\pi^2}n-\vert M(n)\vert \leq A_1n^\frac{1}{2}
\end{equation}

Comparing (\ref{eq17}) and (\ref{eq13}), one can get,

\begin{equation}\label{eq18}
\vert M(n)\vert \leq (A-A_1)n^\frac{1}{2}=Cn^\frac{1}{2}
\end{equation}

where $C$ is constant for a large $n$.

\ 

{\bf Fact 2}.  According to Hardy and Wright (2008),  among the squarefree numbers those for which $\mu(n) = 1$ and those for which $\mu(n) = -1$ occur with about the same frequency for $M(n)=o(n)$. Since $\sum\vert\mu(n)\vert=\frac{6}{\pi^2}n+O(n^{1/2})$, we have,

\begin{equation}\label{ou}
\sum\vert\mu(n)\vert+\sum\mu(n)=2\sum^{n_{\mbox{\tiny e}}} 1=2[\frac{3}{\pi^2}{n_{\mbox{\tiny e}}}+O({n_{\mbox{\tiny e}}}^{1/2})]
\end{equation}

 \begin{equation}\label{ji}
\sum\vert\mu(n)\vert-\sum\mu(n)=-(2\sum^{n_{\mbox{\tiny o}}}-1)=2[\frac{3}{\pi^2}{n_{\mbox{\tiny o}}}+O({n_{\mbox{\tiny o}}}^{1/2})]
\end{equation}

where $n_{\mbox{\tiny e}}$ is squarefree with even number of distinct prime factors, $n_{\mbox{\tiny o}}$  with odd number of distinct prime factors.

Thus, we have,

 \begin{equation}\label{ou2}
\sum\vert\mu(n)\vert+\sum\mu(n)-\frac{6}{\pi^2}n_{\mbox{\tiny e}}\leq A_1n_{\mbox{\tiny e}}^{1/2}
\end{equation}

 \begin{equation}\label{ji2}
\sum\vert\mu(n)\vert-\sum\mu(n)-\frac{6}{\pi^2}n_{\mbox{\tiny o}}\leq A_2n_{\mbox{\tiny o}}^{1/2}
\end{equation}

where $0<A_2< A_1$ are two different constants for a large $n$.

Because those for which $\mu(n) = 1$ and those for which $\mu(n) = -1$ occur with about the same frequency among the squarefree numbers, $n_{\mbox{\tiny e}}=n_{\mbox{\tiny o}}=n'$ for a large $n$.  And from (\ref{ou2}) and (\ref{ji2}), we have,

 \begin{equation}
\sum\mu(n)\leq \frac{A_1-A_2}{2}n'^{1/2}= Cn'^{1/2} \leq Cn^{1/2} 
\end{equation}

where $C$ is a constant for a large $n$.

\ 

{\bf Fact 3}.  From (\ref{eq11}) and (\ref{eq12}),  (\ref{ou}) and (\ref{ji}), when $n$ is large, the following should be true: 

\begin{equation}\label{eq17_1}
\frac{\sum\vert\mu(n)\vert}{n} - \frac{\vert M(n)\vert}{n} = \frac{6}{\pi^2}
\end{equation}

 \begin{equation}\label{ou2_1}
\frac{\sum\vert\mu(n)\vert}{n} - \frac{\sum\mu(n)}{n}=\frac{6}{\pi^2}
\end{equation}

The variations of $\frac{\sum\vert\mu(n)\vert}{n} - \frac{\vert M(n)\vert}{n}$ and $\frac{\sum\vert\mu(n)\vert}{n} - \frac{\sum\mu(n)}{n}$ with $n$ from $n=1$ to $n=2\times 10^7$ are shown in Figure \ref{fig7}, respectively. It can be found, with increasing $n$, the values above are close to the constant $\frac{6}{\pi^2}$.  The absolute error for $\frac{\sum\vert\mu(2\times 10^7)\vert}{2\times 10^7} - \frac{\vert M(2\times 10^7)\vert}{2\times 10^7}$ and $\frac{\sum\vert\mu(2\times 10^7)\vert}{2\times 10^7} - \frac{\sum\mu(2\times 10^7)}{2\times 10^7}$ to $\frac{6}{\pi^2}$ are $-4.6002\times 10^{-5}$ and $4.928\times 10^{-5}$ respectively. 

These above means that the order of the Mertens function $M(n)$ should be about $O(n^\frac{1}{2})$, even $O(n^{\frac{1}{2}+\varepsilon})$ with $0\le \varepsilon <\frac{1}{2}$.  If so, a better result in Ramar$\acute{e}$ (2013) can infer $M(n)$ should be about $O(0.5n^\frac{1}{2})$ for $n\ge10$, even $O(0.1333n^\frac{1}{2})$ for $n\ge 1664$.

Based on these three facts, we conjecture that $M(n)$ should have a similar formula to (\ref{eq10}) or (\ref{eq10_1}) for the upper bound of the $M(n)$ like the following,

\ 

$\exists C>0$ and $C$ is very large, $\exists n_0$, $\forall n>n_0$

\begin{equation}\label{eq_f}
M(n) \leq Cn^{1/2}
\end{equation}

\section{Conclusions}

Based on the results and discussion above, some conclusions can be drawn as follows,

(1) Two elementary formulae for the Mertens function $M(n)$ are obtained, based on the definite recursive relation for M$\ddot{\mathrm{o}}$bius function introduced in Wei (2016). With these formulae, $M(n)$ can be
calculated directly and simply. The most complex operation is only the Mod. However, in the calculation both 
the efficient and/or optimization algorithm for this relation are required when $n$ is large. 

(2) With this relation, $M(1) \sim M(2 \times 10^7) $ are calculated one by one.  Numerical results show that Mertens function $M(n)$ have complicated properties. The sequence of $M(1) \sim M(2 \times 10^7) $ has 19 empirical modes, 16479 zeros, 10043 local maximums/minimums between two neighbor zeros, a maximum of 1240, and a minimum of -1447. 

(3) We also calculated the variation of $\frac{\sum\vert\mu(n)\vert}{n} -\frac{\vert M(n)\vert}{n}$ and $\frac{\sum\vert\mu(n)\vert}{n}-\frac{M(n)}{n}$ from $n=1$ to $n=2\times 10^7$, respectively. Numerical results show that these values are close to the constant $\frac{6}{\pi^2}$ with increasing $n$.

\vspace{8em}

{\Large\bf  Acknowledges}

We thank Dana Jacobsen very much for pointing out the incorrect values about Mertens function $M(n)$ in this e-print and very good comments on the calculating for $M(n)$. We also thank Alisa Sedunova for his very good suggestions.

\vspace{10em}

%
%

\def\thebibliography#1{
{\Large\bf  References}\list
 {}{\setlength\labelwidth{1.4em}\leftmargin\labelwidth
 \setlength\parsep{0pt}\setlength\itemsep{.3\baselineskip}
 \setlength{\itemindent}{-\leftmargin}
 \usecounter{enumi}}
 \def\newblock{\hskip .11em plus .33em minus -.07em}
 \sloppy
 \sfcode`\.=1000\relax}
\let\endthebibliography=\endlist

\end{document}